\documentclass[a4paper,11pt]{article}         
\usepackage{amsmath,amssymb,amsfonts}    
\usepackage{graphicx,caption}            
\usepackage{color}                       
\usepackage{hyperref}                    

\usepackage{makeidx}                  
\usepackage{acronym}
\usepackage{mathrsfs}
\usepackage{float}
\usepackage{vruler}                     
\usepackage{palatino, url, multicol}
\usepackage{setspace}                    
\usepackage{txfonts}

\usepackage{tikz}
\usetikzlibrary{shapes.geometric, arrows}
\usetikzlibrary{calc}

\usepackage[numbers,sort&compress]{natbib}   

\usepackage{import}

\usepackage{caption}

\usepackage{listings}

\usepackage{comment}

\restylefloat{figure}

\begin{document}

\title{The Largest Empty Sphere Problem in 3D Hollowed Point Clouds}
\author{\Large{Netzer Moriya}}

\date{}

\maketitle

\section{Abstract}

We introduce a new approach for the adaptation of the Maximal Internal Envelope method, extended to address the Largest Empty 
Sphere problem within unstructured 3D point clouds. We explore the identification of the Largest Empty Sphere by computing Convex Hull 
vertices and employing a Voidness Score based on Minimal Distance Scoring for optimal segment selection. The integration of Delaunay 
triangulation and Voronoi diagrams facilitates the initial identification of potential Largest Empty Sphere candidates. 
Our analysis reveals the method's efficacy and efficiency, often locating the Largest Empty Sphere in initial computational stages, 
suggesting a lower complexity than initially projected.

\section{Introduction}
The general Largest Empty Sphere Problem (LESP), addresses the problem of finding a hypersphere of largest radius in d-dimensional space 
whose interior does not overlap with any given obstacles. In the case of a 3D point cloud, the problem becomes finding the largest empty 
sphere within the point cloud.
This problem has broad applications across various disciplines, including 3D facility location, 
molecular modeling, and robotics in three-dimensional environments.

Historically, the LESP has been approached from various perspectives, each contributing unique insights into the problem's complexity. 
Aistleitner \cite{Aistleitner:2017} explored the concept in higher dimensions, examining the largest empty box and its volume behavior 
in relation to dimensionality (see also \cite{Chazelle:1986}). In contrast, Shi's work \cite{Shi:1996} focused on a convex optimization 
approach for identifying the largest empty sphere. 
Toussaint \cite{Toussaint:1983} proposed an algorithm for the largest empty circle, integrating locational 
constraints, while Ryshkov \cite{Ryshkov:1988} discussed the empty sphere in the context of lattice coverings, adding a historical 
perspective to the problem.

A common solution technique involves generating a mesh of evenly distributed points throughout the space, then applying a nearest 
neighbors algorithm to determine the furthest mesh point from any given point \( P_i \in P \), thereby identifying the center of the 
largest void. However, this approach can be exteremly computationally demanding in higher dimensions.
An alternative strategy involves Delaunay triangulation, assessing the volume of each formed tetrahedron. These tetrahedrons represent 
potential voids, with the largest among them potentially constituting the largest empty sphere.
A different approach employs the Voronoi diagram on a sphere, combinatorially akin to a 3D convex hull. In this model, each hull face 
correlates to a Voronoi vertex, linked to an empty circle, leading to an efficient \( O(n \log n) \) algorithm \cite{Na:2002}.

The LES problem in shell-like or hollowed structures presents a more intricate challenge, demanding specialized algorithms attuned to their 
unique spatial configurations \cite{Nag:2018}. 
Conventional methods such as Voronoi diagrams and Delaunay triangulation by themselves might not be suffice in such cases. 
Instead, stochastic search approaches \cite{HUGHES:2019} and evolutionary algorithms have proven effective in handling multidimensional spaces, 
offering advantages in execution time and memory usage \cite{Lee:2004}. 
Within the context of metallic hollow spheres, extensive research has been conducted on uniaxial deformation 
behavior under significant deformations, with computational models developed for these scenarios \cite{Speich:2009}. 
Additionally, finite elastostatics within pressurized hollow spheres for compressible materials have been thoroughly investigated, 
yielding explicit closed-form solutions for deformation and stress fields \cite{Abeyaratne:1984}. Further research has provided 
closed-form solutions for hollowed sphere models with specific materials under isotropic loadings, offering valuable reference 
solutions and critical evaluations of macroscopic criteria \cite{Thore:2009}.

This paper aims to extend these methodologies to identify the largest empty sphere within shell-like or hollowed 3D point clouds. 
We introduce a semi-heuristic approach, commencing with Delaunay-Voronoi based void detection, followed by convex hull construction 
for guiding sphere selection. Our methodology ensures no cloud members are enclosed within the identified sphere, addressing the 
unique challenges posed by these structures.

\section{Problem Statement}
Given a set \( P \) of \( n \) points representing a hollowed structure in the Euclidean space, 
we seek to identify the largest sphere within the hollow region that does not enclose any point in $P$. Unlike traditional LESP 
formulations, the hollowed structure presents a unique challenge as the largest empty sphere is expected to be centered within the 
void of the structure, rather than merely avoiding the outer points of $P$. 
Formally, the problem seeks to maximize the radius $r$ of a sphere with center $c$, located within the hollow 
region, ensuring that for all points $p \in P$, the distance from $c$ to $p$ is at least $r$.

\subsection{Mathematical Formulation}
The problem is formalized in the Euclidean space \( \mathbb{R}^3 \). Considering a finite set of points \( P = \{p_1, p_2, \ldots, p_n\} \) 
that outlines a hollowed structure, the aim is to determine the center \( c \) and radius \( r \) of a sphere \( S \in S^{LES}\),
where \( S^{LES} \) denotes the set of potential largest empty spheres in the 3D point cloud, satisfying the following conditions:

\begin{itemize}
    \item The center \( c \) of the sphere \( S \) lies within the void of the structure formed by \( P \), adhering to the spatial 
	constraints imposed by the hollowed configuration.
    \item For every \( p \in P \), the condition \( \| c - p \| \geq r \) holds, where \( \| \cdot \| \) denotes the Euclidean norm. 
	This ensures that \( S \) is external to all points in \( P \).
    \item The objective is to maximize \( r \), subject to the aforementioned conditions, representing the largest feasible sphere 
	within the hollow.
\end{itemize}

Additionally, we introduce an auxiliary function \( \phi: \mathbb{R}^3 \times \mathbb{R} \to \mathbb{R} \) defined 
as \( \phi(c, r) = \min_{p \in P} \| c - p \| - r \). 
The function \( \phi \) quantifies the minimum distance from the sphere boundary to any point in \( P \), subtracting the 
radius \( r \). The LES problem then becomes finding \( (c^*, r^*) = \arg\max_{(c, r)} \phi(c, r) \), subject to \( \phi(c, r) \geq 0 \) 
and other geometric constraints imposed by the hollow structure. To facilitate this, the Minimal Distance Score (MDS) is employed 
to assign a 'Void Score' to various regions within the 3D point cloud. This score aids in pinpointing areas that are in closer proximity 
to the internal cavity of the cloud.

\section{Proposed Methodology}
We propose a new approach to identify the largest empty sphere (LES) within a 3D point cloud, an extension of the 
Maximal Internal Envelope (MIE) method outlined in \cite{Moriya:2023b} from 2D to 3D contexts. This adaptation addresses 
the unique challenges posed by 3D structures, particularly when analyzing voids or cavities.

Initially, the process involves computing the Convex Hull vertices \( P^{CH}_i \) to encapsulate the point cloud. 
The subsequent phase includes a systematic selection of optimal line segments \( L^{best}_{ij} \), which are 
strategically chosen based on their likelihood of intersecting with internal void regions by assigning a 'Voidness Score' to each segment
and choosing the \( L^{best}_{ij} \). 
The approach then integrates Delaunay triangulation with Voronoi diagrams to identify initial points (referred 
to as DV points) for constructing the sets of both MIE points and largest empty spheres. 

This process is iteratively developed using a pseudo-recursive method, enhancing higher-order MIE points detection and potential LES.

Each iteration refines the set of MIE points, facilitating the identification of additional candidate points. 
Starting with a hypothetical sphere centered at a DV point on segment \( L^{best}_{ij} \), this sphere incrementally expands until it 
first intersects with points in the cloud. The sphere's center then moves along \( L^{best}_{ij} \), incrementally expanding in 
predefined step sizes (controlled by a parameter \( 3 \leq k \leq n \)). Points first intersected by the expanding sphere are aggregated 
into set \( S^{LES} \) and recorded as candidate LES spheres, alongside the sphere centers and corresponding radii. This procedure is 
methodically repeated for all segments \( L_{ij} \) associated with \( P^{CH}_i \), thereby ensuring thorough exploration and 
inclusion of significant spheres within the void's interior. 
As the simulation converges to its final stage, the largest sphere, \( \max_{\text{sphere} \in S^{LES}} \text{radius}(\text{sphere}) \), is
then selected as the LES in the 3D point cloud.

We applied this method to two distinct point set configurations: a single sphere with random dispersion along its circumference, 
and a pair of congruent spheres with a deliberately engineered internal void. In both cases, our method successfully revealed the 
internal set of potential largest empty spheres, identifying the sphere with the maximal radius as the LES.

\subsection{A Short Description of the MIE method}

Here we relate to the principle of the method, described in details in \cite{Moriya:2023b} for the 2D case. In a following section, we 
extend the method to the 3D scenario and describe the inclusion of the LES identification as mentioned above.

The proposed methodology for analyzing voids in a three-dimensional point cloud employs a two-stage process, integrating geometric and 
computational techniques.

\textbf{Stage 1: Voids Localization} begins with convex hull analysis to identify principal points of the cloud. 
Using Delaunay Triangulation, Voronoi Diagrams, and the Minimal Distance Scoring (MDS) technique, segments likely to intersect the 
void are determined. We then identify first-order Maximal Internal Envelope (MIE) points based on proximity and angular relations 
to these segments.

\textbf{Stage 2: Iterative Refinement} builds upon the MIE points, enhancing the structural definition of the void iteratively. 
This computationally intensive approach yields a detailed representation of the void's geometry.

\subsubsection{Key Steps and Equations}

\textbf{Minimal Distance Scoring (MDS)} is central to our methodology. The MDS of a segment \(L_{ij}\) is given by:

\begin{equation}
    \text{MDS}(L_{ij}) = \frac{1}{|P \setminus P_{CH}|} \sum_{p_k \in P \setminus P_{CH}} d(L_{ij}, p_k)
\end{equation}

Here, \( d(L_{ij}, p_k) \) represents the distance from a line segment \( L_{ij} \) to a point \( p_k \) in the point set \( P \) 
excluding the points in the Convex Hull \( P_{CH} \). 
This scoring mechanism is fundamental in identifying segments within the point cloud that are likely candidates for forming the 
boundary of the largest void.

In conjunction with MDS, we introduce an auxiliary function \( \phi: \mathbb{R}^3 \times \mathbb{R} \to \mathbb{R} \), defined 
as \( \phi(c, r) = \min_{p \in P} \| c - p \| - r \). 
This function \( \phi \) quantifies the minimum distance from the boundary of a sphere with center \( c \) and radius \( r \) to any 
point in \( P \), ensuring the sphere does not enclose any point in \( P \). 
Therefore, the LES problem is reformulated as an optimization problem where we seek to maximize \( r \) while 
ensuring \( \phi(c, r) \geq 0 \). The MDS aids in this process by providing a 'Voidness Score' to segments, helping in the selection 
of the optimal center \( c \) and radius \( r \) for the largest empty sphere, while adhering to the geometric constraints of the 
hollow structure.

\newpage

\textbf{Stages in Voids Identification} include:
\begin{enumerate}
    \item Convex Hull Identification.
    \item Segment Analysis for intersecting the void using MDS.
    \item Delaunay Triangulation and Voronoi Diagrams to locate the optimal point on the best segment.
\end{enumerate}

\textbf{Void Shape Construction} uses an expanding polygon initiated from the best segment identified by the MDS algorithm. This involves:
\begin{enumerate}
    \item Selecting vertices based on proximity to \( L^{best}_{ij} \).
    \item Iterative expansion and adaptation of the polygon's shape influenced by points in \( P \).
\end{enumerate}

\subsubsection{Complexity Analysis}

The complexity analysis includes:
\begin{itemize}
    \item Convex Hull Calculation: \( O(n \log n) \).
    \item Pair Generation of Convex Hull Vertices: \( O\left(\frac{n^2}{100}\right) \).
    \item Overall Complexity: \( O(n \log n) + 4 \times O\left(\frac{n^2}{100}\right) \).
\end{itemize}

This methodology effectively delineates the shape characteristics of voids in point clouds, balancing detailed geometric analysis 
with computational efficiency.

\subsection{Identification of the Largest Empty Sphere in the Shell-Like 3D Point Cloud}
Given the above methodology for MIE points identification, the expanding spheres intersecting with the MIE points, are aggregated 
into set \( S^{LES} \) and recorded as candidate LES spheres, alongside the sphere centers and corresponding radii. This procedure is 
methodically repeated for all segments \( L_{ij} \) associated with \( P^{CH}_i \), thereby ensuring thorough exploration and 
inclusion of significant spheres within the void's interior. Further orders of the calculations, based on higher order segments 
\cite{Moriya:2023b} are also considered and added to the set \( S^{LES} \).
As the simulation converges to its final stage, the largest sphere, \( \max_{\text{sphere} \in S^{LES}} \text{radius}(\text{sphere}) \), is
then selected as the LES in the 3D point cloud.

\begin{figure}[H]
\begin{center}
\includegraphics[width=7cm]{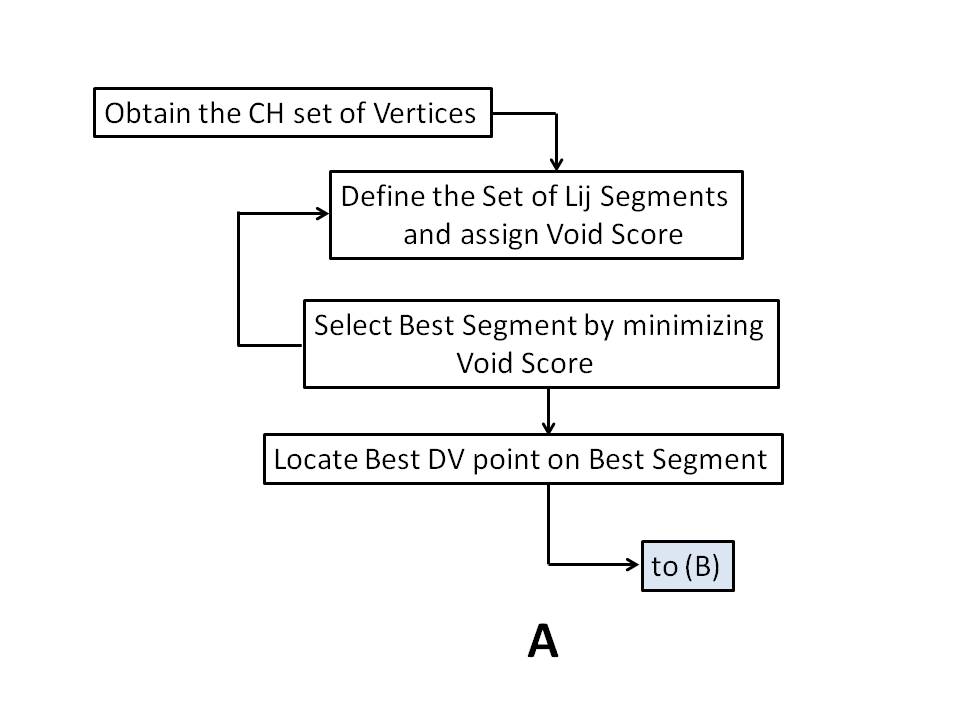} 
\includegraphics[width=7cm]{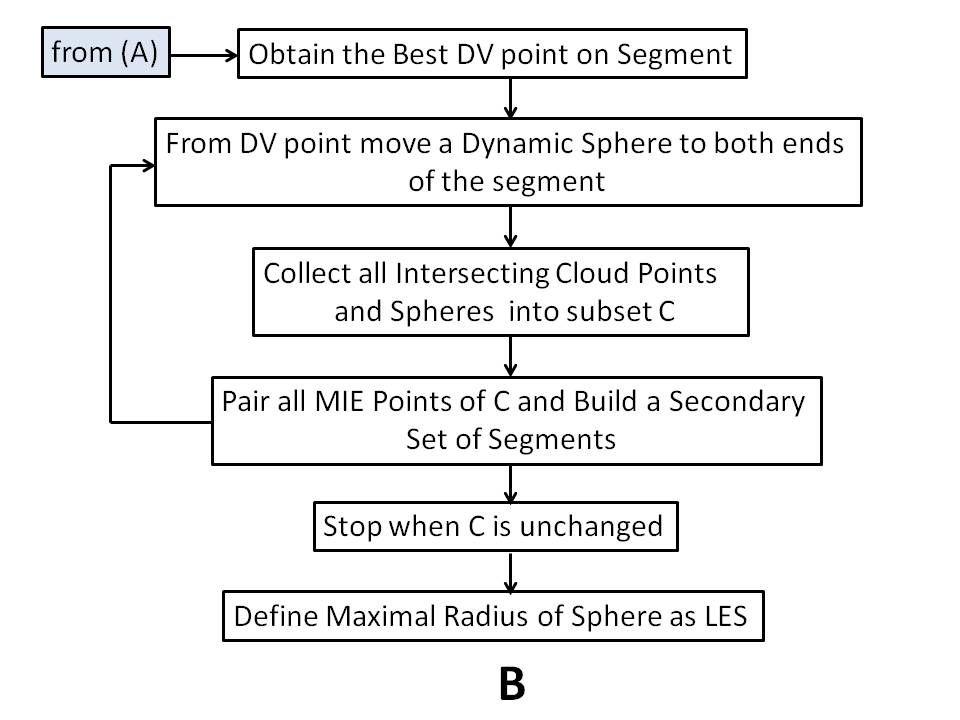}
\caption{Schematic Algorithm for the LES estimation.}
\label{Algorithm_1 Figure}
\end{center}
\end{figure}

The following schematic algorithm \ref{Algorithm_1 Figure} detailing the computational steps for void identification and the 
estimation of the LES in the cloud.

\begin{table}[H]
\centering
\caption{Algorithm for the LES estimation}
\begin{tabular}{rl}
\textbf{Step} & \textbf{Action} \\
1: & Identify Convex Hull \(P^{CH}\) of \(P\) \\
2: & Select a starting point \(P^{CH}_i\) from \(P^{CH}\) \\
3: & For each \(P^{CH}_i\), create and analyze segments \(L_{ij}\) \\
4: & Calculate MDS(\(L_{ij}\)) for each segment \\
5: & Select \(L^{best}_{ij}\) with minimum MDS value \\
6: & Apply Delaunay Triangulation and Voronoi Diagrams to select DV point on \(L^{best}_{ij}\) \\
7: & Along each segment \(L^{best}_{ij}\), define a dynamic Hypothetical Sphere that expands from the DV point. \\
8: & Identify Maximal Internal Envelope (MIE) points from \(P\) that intersect with the expanding circle. \\
9: & Define the Hypothetical Sphere as a LES candidate \\
10: & Incorporate the Hypothetical Sphere into the subset \(S^{LES}\). \\
11: & Update the set \(S^{LES}\) and repeat the process for the new segments formed. \\
12: & Terminate the iterative process when the addition of new MIE points ceases. \\
13: & Select \( \max_{\text{sphere} \in S^{LES}} \text{radius}(\text{sphere}) \) as the LES in the cloud.
\end{tabular}
\label{table:algorithm 1}
\end{table}

\vspace{0.5cm}

\section{Results - Largest Empty Sphere Selection}
In addressing the challenge of identifying the Largest Empty Sphere (LES) within an unstructured 3D point cloud, we have adapted 
and applied the Maximal Internal Envelope (MIE) method. This adaptation involved extending the original MIE approach, as detailed 
in \cite{Moriya:2023b}, to incorporate a third dimension. The process entailed constructing a collection of potential candidates 
for the largest empty sphere. Among these candidates, the sphere that most accurately represents the LES was then determined and 
selected as the optimal solution.

\subsection{Implementation and Analytical Insights}
\begin{figure}[H]
\begin{center}
\includegraphics[width=7cm]{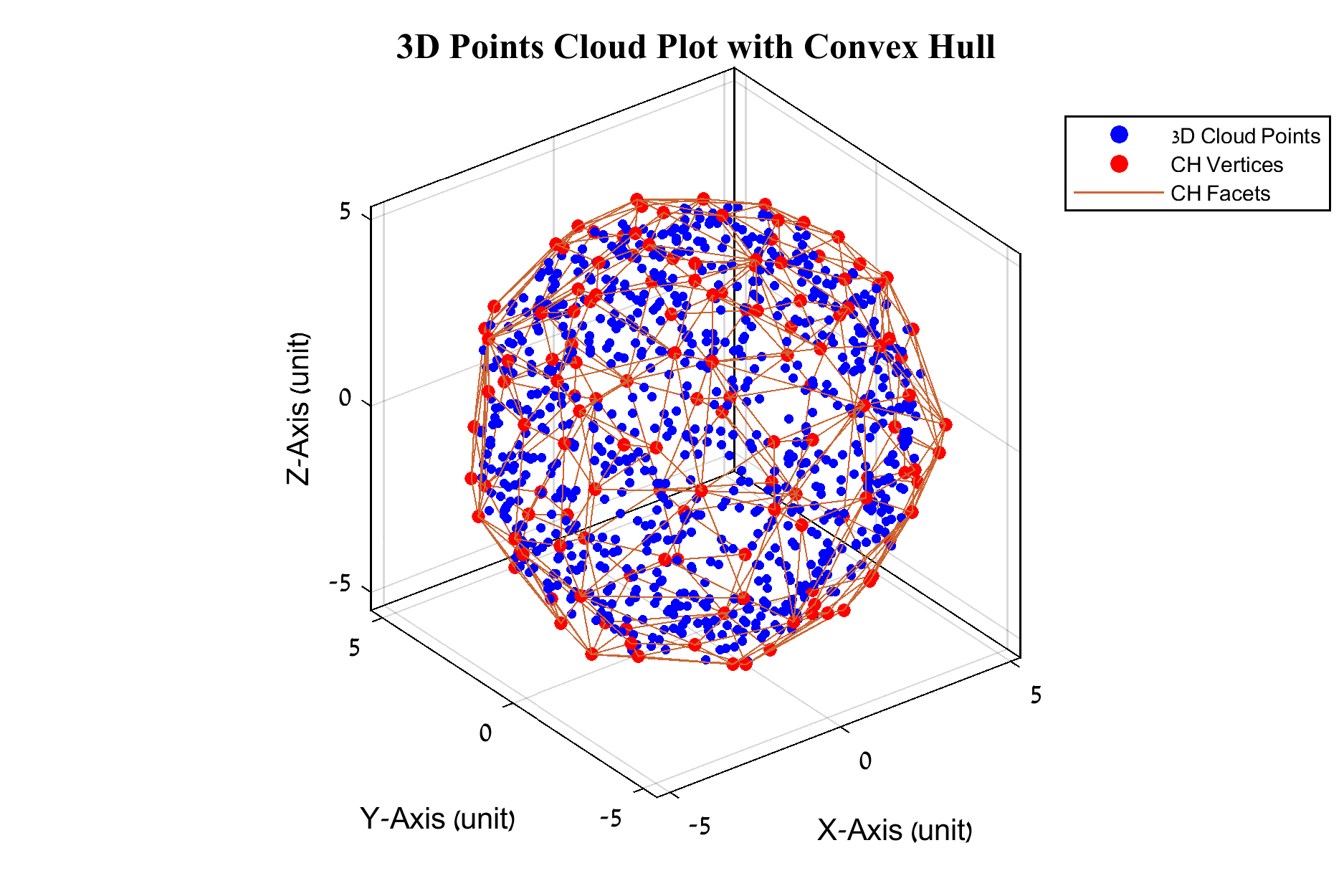}
\includegraphics[width=7cm]{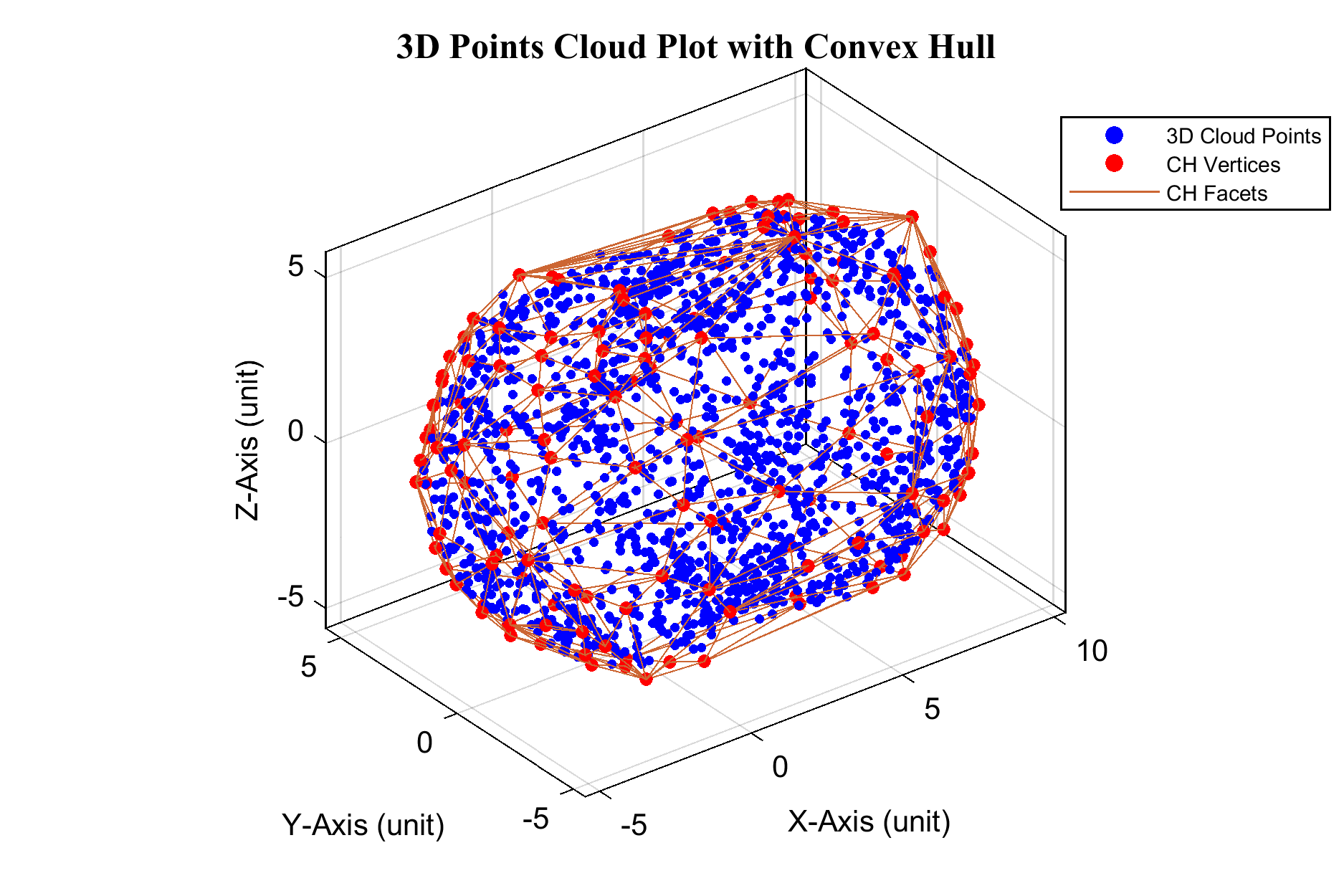}
\caption{The 3D point distributions of a shell-like sphere (left) and the two spheres construction with an internal cavity (right).}
\label{Cloud and CH figure}
\end{center}
\end{figure}

We applied the adapted MIE method to two distinct 3D point cloud configurations: a shell-like sphere and a pair of offset congruent 
spheres. The Convex Hull vertices (\(P^{CH}_i\)) for each configuration were computed, revealing the geometric structure of the point 
clouds (Figure \ref{Cloud and CH figure}).

\begin{figure}[H]
\begin{center}
\includegraphics[width=7cm]{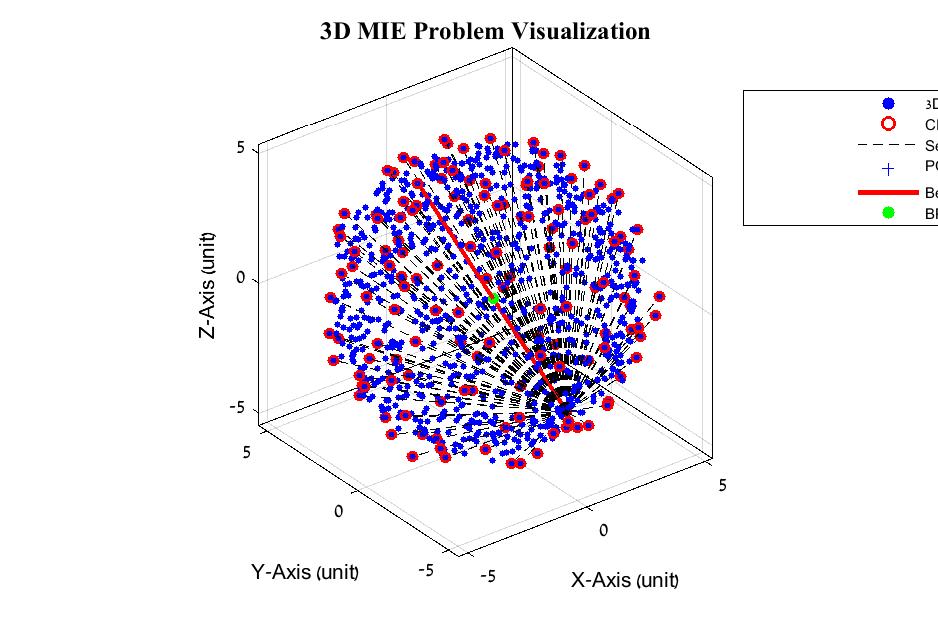}
\includegraphics[width=7cm]{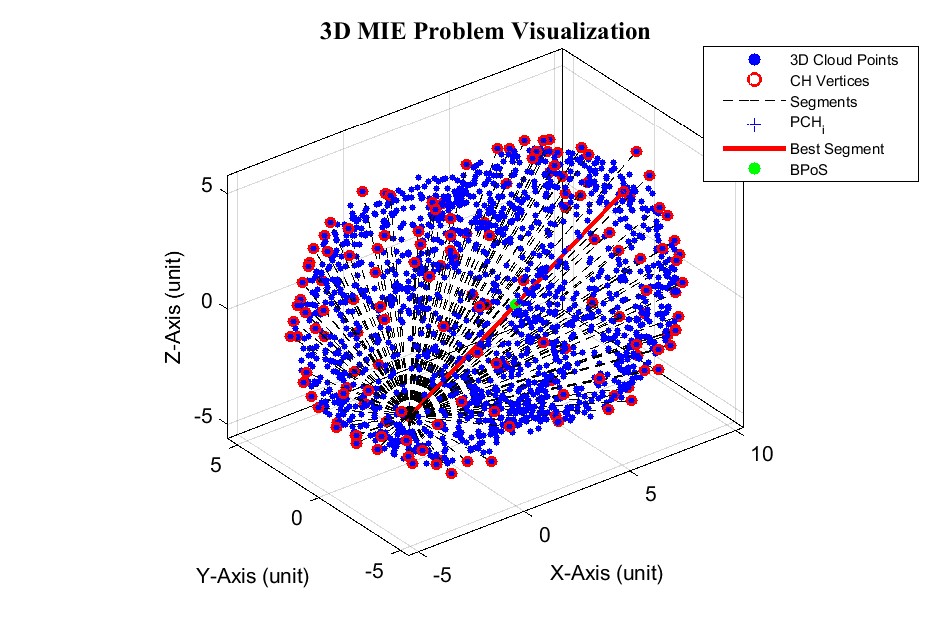}
\caption{Example of best segments selection based on MDS for the two clouds. Every segment is connecting CH pairs of the respective 3D 
points distribution.}
\label{Best segments selection figure}
\end{center}
\end{figure}

\begin{figure}[H]
\begin{center}
\includegraphics[width=7cm]{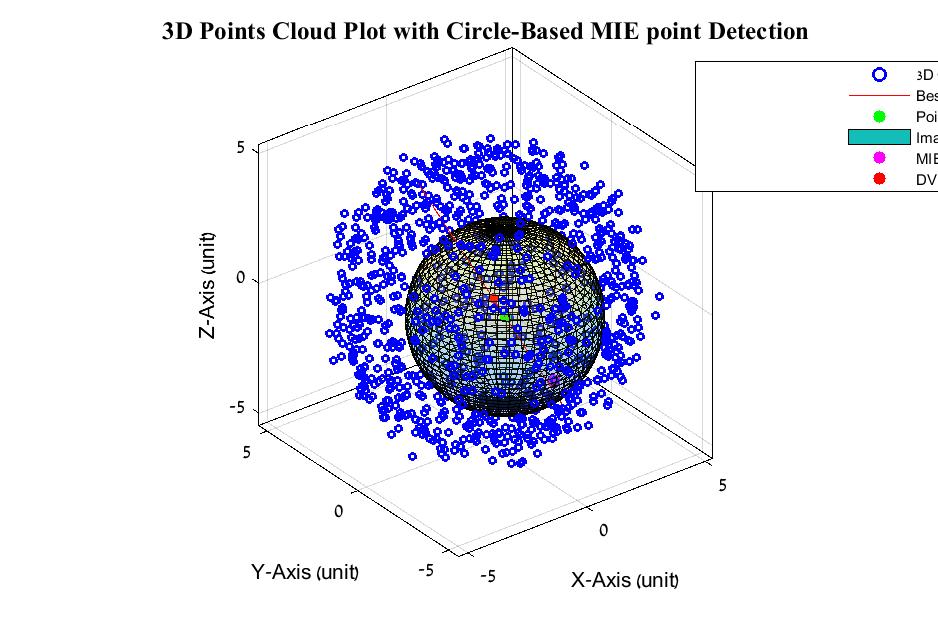}
\includegraphics[width=7cm]{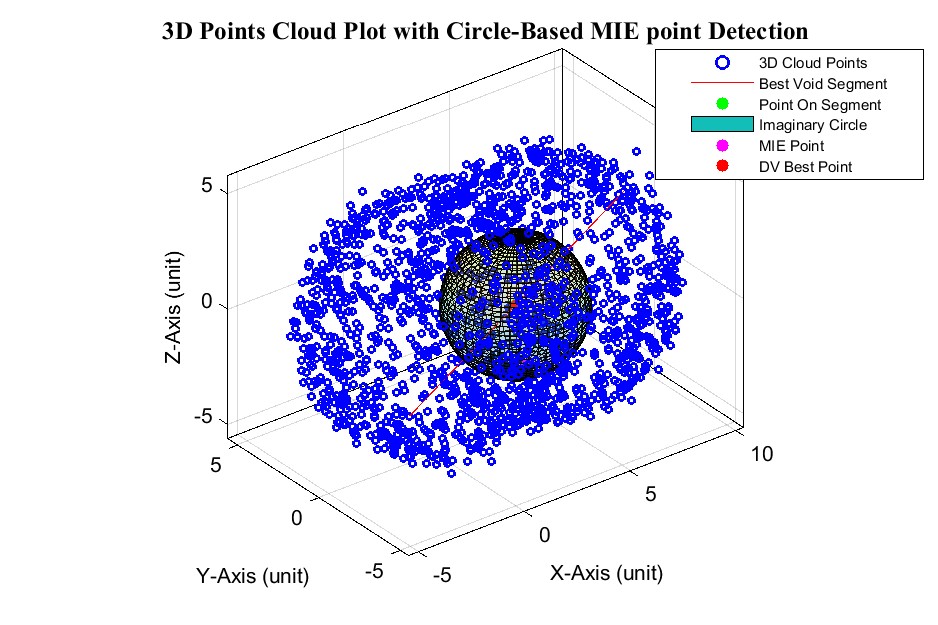}
\caption{Example of MIE points identification based on best segment and VD starting point on segment for the LES selection.}
\label{Imaginary Sphere figure}
\end{center}
\end{figure}

For each configuration, optimal line segments (\(L^{best}_{ij}\)) were identified based on a 'Voidness Score', derived from Minimal 
Distance Scoring (MDS). Delaunay triangulation and Voronoi diagrams facilitated the identification of initial points for LES estimation. 
The dynamic projection of hypothetical spheres from these points and their interaction with the cloud were analyzed to determine the 
potential LES candidates (Figures \ref{Best segments selection figure} and \ref{Imaginary Sphere figure}).

The process involved projecting a hypothetical, dynamic sphere from the Delaunay-Voronoi (DV) point along the selected segment. Points 
intersecting with the radius changing sphere, as shown in Figure \ref{Imaginary Sphere figure}, were included in the subset 
\( C^{MIE} \subseteq P \) with all potential largest empty spheres include in \( S^{LES} \). 
These points, determaning the radii of the potential Largest Empty Spheres, were then subject to maximal radius analysis and upon 
identifying the largest radius, the corresponding sphere was identified as the best estimation for the LES.

\subsection{Maximal Radius Analysis and LES Selection}
\begin{figure}[H]
\begin{center}
\includegraphics[width=7cm]{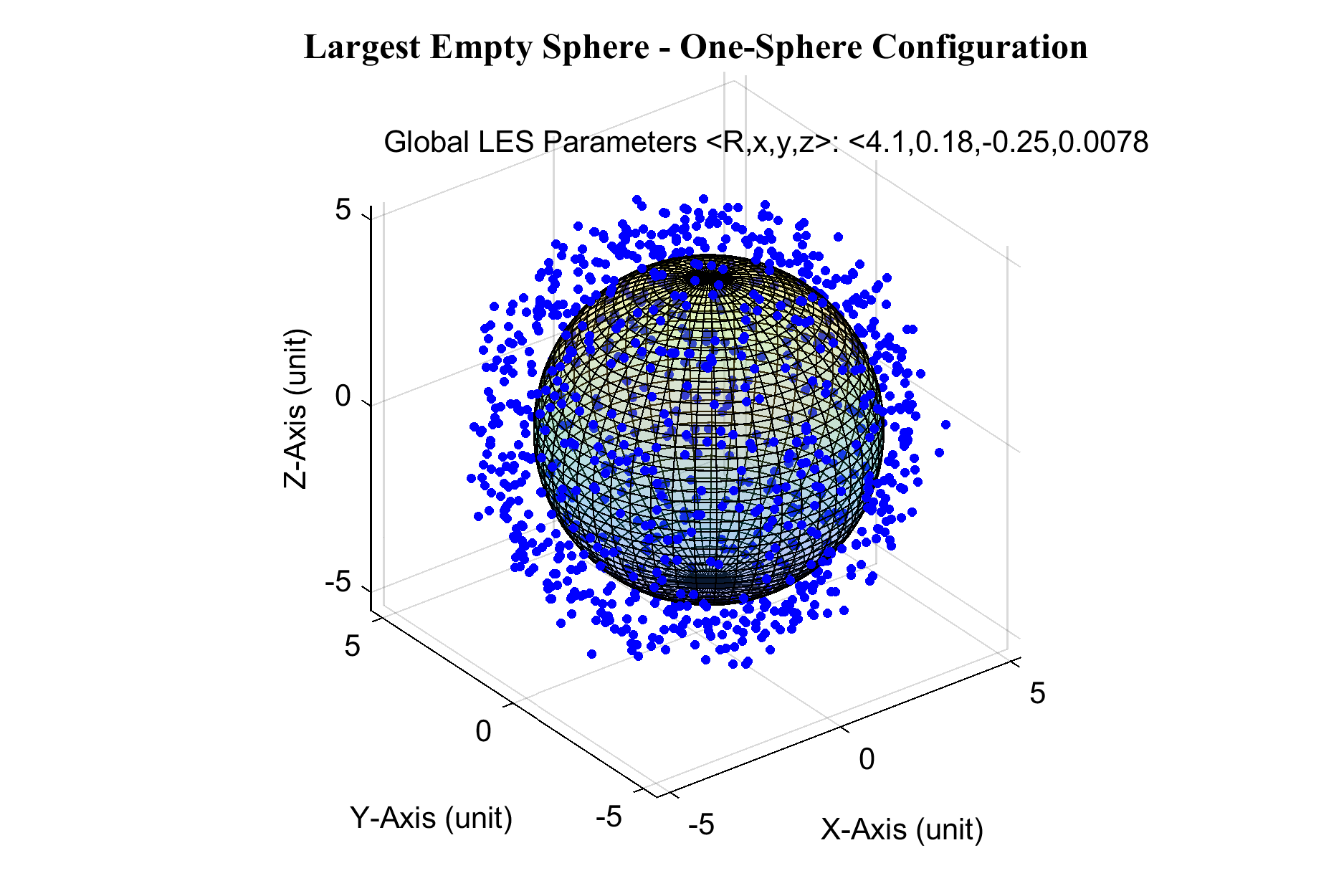}
\includegraphics[width=7cm]{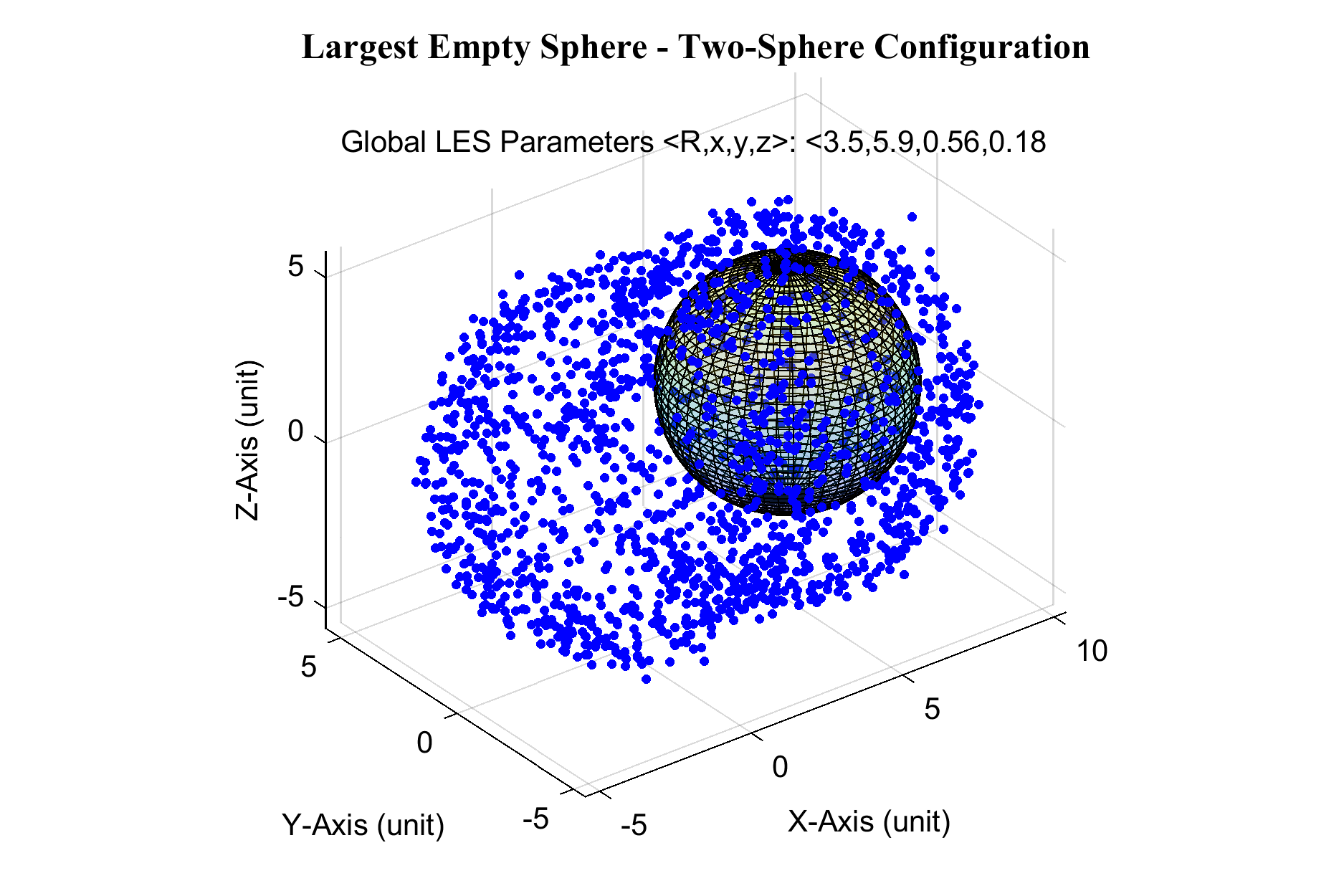}
\caption{The best sphere entitled LES for the shell-like sphere (left) and the two-spheres configuration (right)} 
\label{LES figure}
\end{center}
\end{figure}

The collection of potential LES candidates (\(S^{LES}\)) was subjected to a maximal radius analysis. The largest sphere from \(S^{LES}\) 
was identified as the LES (Figure \ref{LES figure}).

\subsection{Computational Efficiency and Observations}
Our analysis revealed the MIE method's efficacy and efficiency in identifying the LES within diverse 3D point cloud configurations. 
Interestingly, the LES was often identified in the initial computation phase, indicating a lower computational complexity than 
initially projected.

\vspace{0.4cm}

\textbf{Key Insights:}

\begin{itemize}
    \item \textbf{Efficiency and Robustness:} 
	The method of Maximal Internal Envelope has been proved both eficient and robust method to identify that Largest Empty Sphere 
	in cases where cavity is enclosed in a 3D point cloud.
    \item \textbf{Fast Convergence:} 
	In all our observations, over a vast number of simulations obtained with numerious different 3D point cloud configurations, the 
	LES was identified during the first order of the calculationd, i.e., there was no need for higher order iterations as the
	largest empty sphere was almost imediatly identified during the CH stage. 
	\item \textbf{Computational Complexity:} 
	The overall expected complexity analysis is substantially lower than the $O(n \log n) + 4 \times O\left(\frac{n^2}{100}\right)$ 
	deduced for the MIE case. For the LES problem, our assumption is that only the first term should be considered..
\end{itemize}

\section{Conclusions}

The Largest Empty Sphere Problem involving 3D hollowed structures, is an important extension of the classical problem, representing key 
concepts in spatial analysis and optimization of unstructured point clouds in 3D.

This study demonstrates the adaptation of the Maximal Internal Envelope (MIE) methodology, explicitly tailored to tackle the challenge 
of identifying the Largest Empty Sphere (LES) within unstructured 3D point clouds. 
The approach described in this paper involves the computation of Convex Hull vertices, followed by the application of a Voidness Score 
derived from the Minimal Distance Scoring (MDS) algorithm, which is instrumental in the optimal selection of segments. 
Furthermore, the incorporation of Delaunay triangulation and Voronoi diagrams plays a pivotal role in the preliminary identification 
of prospective LES candidates. 
Empirical analysis of this method, proved efficient in identifying the LES within diverse 3D point cloud 
configurations. Notably, the LES was often identified in the initial computational phase, indicating a lower computational complexity than 
previously estimated, prior to the initial simulation experiments. 
Future work may focus on enhancing the algorithm's efficiency and exploring its applicability to larger, more intricate 
datasets in 3D facility location, molecular modeling, and robotics in three-dimensional environments and related fields.

\section*{Declarations}
All data-related information and coding scripts discussed in the results section are available from the 
corresponding author upon request.

\bibliographystyle{unsrtnat}

\renewcommand{\bibname}{References}

\end{document}